\documentclass[a4paper,12pt]{article}

\usepackage[top=2.5cm,bottom=2.5cm,left=2.25cm,right=2.25cm]{geometry}
\usepackage{amsmath,amsthm,mathrsfs,graphicx,amsfonts}
\usepackage{bm}
\usepackage{cases}
\usepackage[english]{babel}
\usepackage{latexsym}
\usepackage{txfonts}
\usepackage[numbers,sort&compress]{natbib}
\usepackage[natural]{xcolor}
\usepackage{rotating}
\usepackage{mathtools}
\usepackage{enumitem}

\makeatletter
\def\@maketitle{
  \newpage
  \null
  \vskip 2em
  \begin{center}
  \let\footnote\thanks
    {\LARGE \@title \par}
    \vskip 1.5em
    {\large
      \lineskip .5em
      \parbox{\textwidth}{\centering\@author}
      \par
    }
    \vskip 1em
    {\large \@date}
  \end{center}
  \par
  \vskip 1.5em}
\makeatother

\newtheorem{lem}{Lemma}[section]
\newtheorem{thm}[lem]{Theorem}

\newtheorem{pro}[lem]{Proposition}

\newtheorem{claim}[lem]{\indent Claim}

\usepackage[colorlinks=false, pdfborder={0 0 0}]{hyperref}

\makeatletter
\@removefromreset{lem}{section}
\makeatother

\title{On the number of edge-disjoint $K_r$ in $K_{r+1}$-free graphs}

\author{Zhen Liu\footnote{Email: 1552580575@qq.com}, ~Qinghou Zeng\footnote{Research supported by National Key R\&D Program of China (Grant No. 2023YFA1010202) and National Natural Science Foundation of China (Grant No. 12371342). Email: zengqh@fzu.edu.cn (Corresponding
author)}\\
{\small Center for Discrete Mathematics, Fuzhou University, Fujian, 350003, China}}
\date{}
\begin{document}
\maketitle

\begin{abstract}
Gy\H ori and Keszegh [\emph{Combinatorica} \textbf{37} (2017), 1113--1124] proved that every $K_4$-free graph on $n$ vertices with $\lfloor n^2/4\rfloor+k$ edges contains at least $k$ pairwise edge-disjoint triangles. We extend this result to cliques of arbitrary fixed order. Specifically, we prove that for every fixed integer $r\ge 3$ and all sufficiently large $n$, every $K_{r+1}$-free graph on $n$ vertices with $t_{r-1}(n)+k$ edges contains at least $k$ pairwise edge-disjoint copies of $K_r$, where $t_{r-1}(n)$ denotes the number of edges in the Tur\'an graph $T_{r-1}(n)$. Moreover, equality holds if and only if $G$ is obtained from $T_{r-1}(n)$ by adding edges within one of its parts, with the added edges forming a triangle-free graph.
\end{abstract}

\section{Introduction}

Throughout, graphs are finite and simple. Let $\nu_r(G)$ denote the
maximum number of pairwise edge-disjoint copies of $K_r$ in a graph $G$.
Write $T_{r-1}(n)$ for the balanced complete $(r-1)$-partite graph on
$n$ vertices and put $t_{r-1}(n)=e(T_{r-1}(n))$. By Tur\'an's theorem,
every graph with more than $t_{r-1}(n)$ edges contains $K_r$.

Erd\H{o}s~\cite{Erdos1971} posed the problem of
determining the largest integer $k$ such that every
graph on $n$ vertices with $t_2(n)+k$ edges contains
$k$ pairwise edge-disjoint triangles.
Gy\H{o}ri~\cite{Gyori1988} solved this problem for
sufficiently large $n$, with the threshold for odd $n$
subsequently corrected in~\cite{Gyori1992}.
The largest such $k$ is $2n-10$ when $n$ is odd
and $3n/2-5$ when $n$ is even.
Gy\H ori~\cite{Gyori1991} extended this result to
larger cliques. He proved that, for every fixed
$r\ge4$ and all sufficiently large $n$, every
graph $G$ on $n$ vertices with $t_{r-1}(n)+k$
edges contains at least $k$ pairwise edge-disjoint
copies of $K_r$, provided that
$k\le 3\left\lfloor\frac{n+1}{r-1}\right\rfloor-5$.
He also proved that, for fixed $r\ge3$ and
$k=o(n^2)$, $\nu_r(G)\ge k-O(k^2/n^2)$ in~\cite{Gyori1992}.
For unrestricted $k$, Gy\H ori and
Tuza~\cite{GyoriTuza1987} established
$\nu_r(G)\ge 2k/[r(r-2)]$ for $r\ge4$.
For $r=3$, Kr\'al', Lidick\'y, Martins and
Pehova~\cite{KralLidickyMartinsPehova2019}
established the asymptotic bound
$\nu_3(G)\ge 2k/3-o(n^2)$.
Subsequently, Blumenthal
et al.~\cite{BlumenthalEtAl2021}
determined the exact extremal value of a
related edge-decomposition parameter,
which implies that
$\nu_3(G)\ge(2k-1)/3$ for all sufficiently
large $n$.
More recently, Balogh and
Wigal~\cite{BaloghWigal2025} resolved
a conjecture of Gy\H ori by showing that,
for every fixed $r\ge3$,
$
\nu_r(G)\ge\left(\frac{2}{r}-o(1)\right)k.
$
The leading constant $2/r$ is asymptotically
best possible.

For $K_4$-free graphs, Gy\H ori~\cite{Gyori1992}
conjectured that every graph on $n$ vertices
with $t_2(n)+k$ edges contains at least $k$
pairwise edge-disjoint triangles, without any
restriction on $k$. Huang and
Shi~\cite{HuangShi2014} made progress towards
this conjecture by proving that
$\nu_3(G)\ge32k/35$.
The conjecture was eventually confirmed by
Gy\H ori and Keszegh.

\begin{thm}[Gy\H ori and Keszegh~\cite{GyoriKeszegh2017}]
\label{thm:gyori-keszegh}
Every $K_4$-free graph $G$ on $n$ vertices
with $\lfloor n^2/4\rfloor+k$ edges contains
at least $k$ pairwise edge-disjoint triangles.
\end{thm}

Recently, Zhang~\cite{Zhang2026} characterized
the graphs attaining equality in
Theorem~\ref{thm:gyori-keszegh} for
sufficiently large $n$.
In this paper, we extend
Theorem~\ref{thm:gyori-keszegh} to
$K_{r+1}$-free graphs for every fixed $r\ge3$
and characterize all graphs attaining equality.
Let $\mathcal E_r(n)$ denote the family of
graphs obtained from $T_{r-1}(n)$ by adding
edges within one of its parts, with the added
edges forming a triangle-free graph. The set
of added edges may be empty.

\begin{thm}\label{thm:main}
For every fixed integer $r\ge3$ and all
sufficiently large $n$, every $K_{r+1}$-free
graph $G$ on $n$ vertices with $t_{r-1}(n)+k$
edges contains at least $k$ pairwise
edge-disjoint copies of $K_r$.
Moreover, equality holds if and only if
$G\in\mathcal E_r(n)$.
\end{thm}

\noindent\textbf{Remark 1.}
By Tur\'an's theorem, every $n$-vertex
$K_{r+1}$-free graph with $t_{r-1}(n)+k$
edges satisfies
$k\le t_r(n)-t_{r-1}(n)$,
with equality attained by $T_r(n)$.
In contrast, by Mantel's theorem, the possible
numbers of edges in graphs belonging to
$\mathcal E_r(n)$ are precisely
$t_{r-1}(n)+k$, where $k$ is an integer
satisfying
$0\le k\le
\lfloor\lceil n/(r-1)\rceil^2/4\rfloor$.
Consequently, Theorem~\ref{thm:main} implies that,
for sufficiently large $n$, every $n$-vertex
$K_{r+1}$-free graph with $t_{r-1}(n)+k$ edges
contains at least $k+1$ pairwise edge-disjoint
copies of $K_r$ whenever
$\lfloor\lceil n/(r-1)\rceil^2/4\rfloor
<k\le t_r(n)-t_{r-1}(n)=\frac{n^2}{2r(r-1)}+O_r(1)$.

\noindent\textbf{Remark 2.}
For $r=3$, Gy\H ori and
Keszegh~\cite{GyoriKeszegh2017} proved that the
lower bound in Theorem~\ref{thm:main} holds for
every $n$. However, the assumption that $n$ is
sufficiently large cannot be omitted for $r\ge4$.
Indeed, consider $G=K_{2,2,2,2}$, which is
$K_5$-free and satisfies
$e(G)=24=t_3(8)+3$, but $\nu_4(G)=2$.
For every $r\ge5$, we will construct in
Section~\ref{proof} a $K_{r+1}$-free graph
on $n=(r-1)(r-2)+1$ vertices with
$t_{r-1}(n)+2r-6$ edges that contains fewer
than $2r-6$ pairwise edge-disjoint copies
of $K_r$.

\section{Notation and preliminary results}
\label{sec:preliminaries}

Throughout, all graphs are finite and simple, and
$r\ge3$ is fixed unless otherwise specified.
For every positive integer $a$, write
$[a]=\{1,\ldots,a\}$.
For a graph $G$, let $V(G)$ and $E(G)$ denote
its vertex and edge sets, respectively, and write
$e(G)=|E(G)|$.
For a vertex $v\in V(G)$, let $N_G(v)$ denote
its neighborhood, $d_G(v)=|N_G(v)|$ its degree,
and $\Delta(G)$ the maximum degree of $G$.
For $X\subseteq V(G)$, let $G[X]$ denote the
subgraph of $G$ induced by $X$.
A subgraph $H$ of $G$ is called a
\emph{spanning subgraph} if $V(H)=V(G)$.
For $v\in V(G)$ and $X\subseteq V(G)$, write
$d_G(v,X)=|N_G(v)\cap X|$.
For a set $S\subseteq V(G)$, let
$N_G(S)=\bigcup_{v\in S}N_G(v)$ denote the
set of vertices adjacent to at least one
vertex of $S$. In particular,
$N_G(\varnothing)=\varnothing$.
For disjoint sets $X,Y\subseteq V(G)$,
let $e_G(X,Y)$ denote the number of edges
between $X$ and $Y$.
We omit the subscript $G$ when the graph
is clear.

An $(r-1)$-partition of an $n$-vertex graph
$G$ is a partition $\mathcal P=(V_1,\ldots,V_{r-1})$
of $V(G)$ into $r-1$ nonempty, pairwise disjoint
sets such that $\bigcup_{i\in[r-1]}V_i=V(G)$.
The partition $\mathcal P$ is \emph{balanced} if
$\bigl||V_i|-|V_j|\bigr|\le1$ for all
$i,j\in[r-1]$.
With respect to $\mathcal P$, an edge of $G$
is called \emph{internal} if both endpoints
belong to the same part, and \emph{crossing}
if its endpoints belong to distinct parts.
Define

\[
\begin{gathered}
I=I(\mathcal P)=\sum_{i=1}^{r-1}e(G[V_i]),
\qquad
M=M(\mathcal P)=\sum_{1\le i<j\le r-1}
\bigl(|V_i||V_j|-e_G(V_i,V_j)\bigr),\\
D=D(\mathcal P)=t_{r-1}(n)
-\sum_{1\le i<j\le r-1}|V_i||V_j|.
\end{gathered}
\]
Thus $I$ is the number of internal edges,
$M$ is the number of missing crossing edges,
and $D\ge0$, with equality if and only if
$\mathcal P$ is balanced.

For each $i\in[r-1]$ and $v\in V_i$, define
\[
\mu(v)=\sum_{j\in[r-1]\setminus\{i\}}
\bigl(|V_j|-d_G(v,V_j)\bigr).
\]
Thus $\mu(v)$ is the number of nonneighbors
of $v$ outside $V_i$. Since each missing
crossing edge is counted at both endpoints,
we have
\begin{equation}\label{eq:sum-deficiency}
\sum_{v\in V(G)}\mu(v)=2M.
\end{equation}

Suppose that $\mathcal P$ maximizes the
number of crossing edges among all
$(r-1)$-partitions of $G$. Then, for all
distinct $i,j\in[r-1]$ and every $v\in V_i$,
we have
\begin{equation}\label{eq:local-optimal}
d_G(v,V_i)\le d_G(v,V_j).
\end{equation}
Indeed, the inequality is immediate if
$|V_i|=1$. Otherwise, moving $v$ from
$V_i$ to $V_j$ changes the number of
crossing edges by
$d_G(v,V_i)-d_G(v,V_j)$, which is
nonpositive by the maximality of
$\mathcal P$.

Let $\mathcal K_r(G)$ denote the family of all
copies of $K_r$ in $G$. A \emph{fractional
$K_r$-packing} is a function
$w:\mathcal K_r(G)\to[0,\infty)$ satisfying
\[
\sum_{\substack{K\in\mathcal K_r(G)\\e\in E(K)}}
w(K)\le1
\qquad\text{for every }e\in E(G).
\]
The maximum total weight of a fractional
$K_r$-packing is denoted by
$\nu_r^*(G)=\max_w\sum_{K\in\mathcal K_r(G)}w(K)$.
A \emph{fractional $K_r$-cover} is a function
$z:E(G)\to[0,1]$ satisfying
$\sum_{e\in E(K)}z(e)\ge1$ for every
$K\in\mathcal K_r(G)$.
By linear programming duality,
\[
\nu_r^*(G)=\min_z\sum_{e\in E(G)}z(e),
\]
where the minimum is taken over all fractional
$K_r$-covers. The restriction $z(e)\le1$
causes no loss of generality, since any
weight greater than one can be reduced to one.

We shall use the following three results.

\begin{thm}[Haxell--R\"odl~\cite{HaxellRodl2001}]
\label{thm:haxell-rodl}
For every fixed integer $r\ge3$, every graph
$G$ on $n$ vertices satisfies
$\nu_r^*(G)-\nu_r(G)=o(n^2)$.
\end{thm}

\begin{thm}[F\"uredi~\cite{Furedi2015}]
\label{thm:furedi}
For every integer $r\ge3$, every $n$-vertex
$K_r$-free graph $H$ can be made
$(r-1)$-partite by deleting at most
$t_{r-1}(n)-e(H)$ edges.
\end{thm}

We introduce some basic definitions
for hypergraphs. A hypergraph $\mathcal H$ consists
of a finite nonempty vertex set $V(\mathcal H)$
and a set $E(\mathcal H)$ of nonempty subsets
of $V(\mathcal H)$, called hyperedges.
It is \emph{linear} if $|e\cap f|\le1$
for all distinct $e,f\in E(\mathcal H)$.
The \emph{rank} of $\mathcal H$ is
$\max\{|e|:e\in E(\mathcal H)\}$, taken to be
zero if $E(\mathcal H)=\varnothing$.
For each $x\in V(\mathcal H)$, define
$d_{\mathcal H}(x)=
|\{e\in E(\mathcal H):x\in e\}|$
and
$\Delta(\mathcal H)=
\max_{x\in V(\mathcal H)}d_{\mathcal H}(x)$.

An \emph{edge-list assignment} $L$ assigns
a finite set of available colours $L(e)$
to each $e\in E(\mathcal H)$.
A \emph{proper $L$-edge-colouring} is a function
$\varphi:E(\mathcal H)\to
\bigcup_{e\in E(\mathcal H)}L(e)$
such that $\varphi(e)\in L(e)$ for every
$e\in E(\mathcal H)$, and
$\varphi(e)\ne\varphi(f)$ whenever
$e,f\in E(\mathcal H)$ are distinct
and $e\cap f\ne\varnothing$.

We use the following form of Kahn's
list-edge-colouring theorem.
\begin{thm}[Kahn~\cite{Kahn1996}]
\label{thm:kahn}
For every integer $k\ge2$ and every $\gamma>0$,
there exists an integer $D_0=D_0(k,\gamma)$
such that the following holds.
Let $\mathcal H$ be a linear hypergraph of rank
at most $k$ with $\Delta(\mathcal H)\ge D_0$.
If $L$ is an edge-list assignment satisfying
$|L(e)|\ge(1+\gamma)\Delta(\mathcal H)$
for every $e\in E(\mathcal H)$, then
$\mathcal H$ admits a proper $L$-edge-colouring.
\end{thm}

\section{Fractional stability and extension lemmas}
\begin{lem}\label{lem:fractional-stability}
Fix an integer $r\ge3$, and let $G$ be an
$n$-vertex $K_{r+1}$-free graph satisfying
$\nu_r(G)\le e(G)-t_{r-1}(n)+o(n^2)$.
Let $\mathcal P$ be an $(r-1)$-partition of $G$
maximizing the number of crossing edges.
Then $M(\mathcal P)+D(\mathcal P)=o(n^2)$.
\end{lem}

\begin{proof}
Let $z:E(G)\to[0,1]$ be a fractional $K_r$-cover minimizing
$\sum_{e\in E(G)}z(e)$.
Thus $\sum_{e\in E(C)}z(e)\ge1$ for every copy $C$ of
$K_r$ in $G$.
By linear programming duality,
$\sum_{e\in E(G)}z(e)=\nu_r^*(G)$.
Theorem~\ref{thm:haxell-rodl} and the hypothesis imply
\begin{equation}\label{eq:fs-weight-total}
\begin{aligned}
\sum_{e\in E(G)}(1-z(e))=e(G)-\nu_r^*(G)=e(G)-\nu_r(G)-o(n^2)
 \ge t_{r-1}(n)-o(n^2).
\end{aligned}
\end{equation}
\begin{claim}\label{clm:fs-rounding-error}
We have $\sum_{e\in E(G)}z(e)(1-z(e))=o(n^2)$.
\end{claim}

\begin{proof}
Fix any $\varepsilon>0$. By \eqref{eq:fs-weight-total}
and $t_{r-1}(n)=\frac{r-2}{2(r-1)}n^2+O(1)$,
for all sufficiently large $n$,
\begin{equation}\label{eq:fs-small-error}
\frac1{n^2}\sum_{e\in E(G)}(1-z(e))
\ge \frac{r-2}{2(r-1)}-\varepsilon^2.
\end{equation}

We first show that there exists a function
$w:V(G)\to[0,1]$ satisfying
$\sum_{v\in V(G)}w(v)=1$, such that the vertices
with positive weight form a clique and
\begin{equation}\label{eq:fs-weight-reduction}
\begin{aligned}
&\frac1{n^2}\sum_{e\in E(G)}
(1-z(e))(1+\varepsilon z(e))\le
\sum_{uv\in E(G)}
(1-z(uv))(1+\varepsilon z(uv))w(u)w(v).
\end{aligned}
\end{equation}

Initially, let $w(v)=1/n$ for every $v\in V(G)$.
Then $\sum_{v\in V(G)}w(v)=1$, and equality holds in
\eqref{eq:fs-weight-reduction}.
Suppose that $u,v\in V(G)$ are nonadjacent vertices
with $w(u),w(v)>0$. Interchanging $u$ and $v$ if
necessary, we may assume that
\[
\begin{aligned}
&\sum_{a\in N_G(u)}
(1-z(ua))(1+\varepsilon z(ua))w(a)\ge
\sum_{a\in N_G(v)}
(1-z(va))(1+\varepsilon z(va))w(a).
\end{aligned}
\]
Define a new function $w':V(G)\to[0,1]$ by
$w'(u)=w(u)+w(v)$, $w'(v)=0$, and
$w'(a)=w(a)$ for every
$a\in V(G)\setminus\{u,v\}$.
Since $uv\notin E(G)$, we have
\[
\begin{aligned}
&\sum_{ab\in E(G)}
(1-z(ab))(1+\varepsilon z(ab))
\bigl(w'(a)w'(b)-w(a)w(b)\bigr)\\
&=w(v)\left(
\sum_{a\in N_G(u)}
(1-z(ua))(1+\varepsilon z(ua))w(a)
-
\sum_{a\in N_G(v)}
(1-z(va))(1+\varepsilon z(va))w(a)
\right)\ge0.
\end{aligned}
\]
Thus replacing $w$ by $w'$ preserves the total weight
and does not decrease the right-hand side of
\eqref{eq:fs-weight-reduction}. Moreover, since
$w(u),w(v)>0$, this operation reduces the number of
vertices with positive weight by exactly one.
Starting with $n$ vertices of positive weight, we
repeat this operation whenever two nonadjacent vertices
have positive weight. Since each operation reduces the
number of such vertices by one, the process terminates
after at most $n-1$ operations. The resulting function
$w$ still satisfies \eqref{eq:fs-weight-reduction}, and
any two vertices with positive weight are adjacent.
Since $G$ is $K_{r+1}$-free, at most $r$ vertices have
positive weight, as required.

Choose distinct vertices $v_1,\ldots,v_r\in V(G)$
containing all vertices with positive weight, and label
them so that
\[
0\le w(v_1)\le w(v_2)\le\cdots\le w(v_r),
\qquad \sum_{i=1}^r w(v_i)=1.
\]
This is possible for $n\ge r$, since at most $r$
vertices have positive weight.

Define an auxiliary function $z'$ on all unordered pairs
of $\{v_1,\ldots,v_r\}$ by
\[
z'(v_iv_j)=
\begin{cases}
z(v_iv_j), & \text{if }w(v_i)w(v_j)>0,\\
1, & \text{otherwise},
\end{cases}
\qquad 1\le i<j\le r.
\]
The first case is well-defined because the vertices with
positive weight form a clique. Clearly,
$0\le z'(v_iv_j)\le1$. Moreover,
$\sum_{1\le i<j\le r}z'(v_iv_j)\ge1$.
Indeed, if all $r$ weights are positive, this follows
from the fractional $K_r$-cover condition. Otherwise,
at least one value of $z'$ equals $1$.

Since $w$ satisfies \eqref{eq:fs-weight-reduction}
and all vertices outside $\{v_1,\ldots,v_r\}$ have
weight zero, we obtain
\begin{equation}\label{eq:fs-weight-comparison}
\begin{aligned}
&\frac1{n^2}\sum_{e\in E(G)}
(1-z(e))(1+\varepsilon z(e))\le
\sum_{1\le i<j\le r}
(1-z'(v_iv_j))(1+\varepsilon z'(v_iv_j))
w(v_i)w(v_j).
\end{aligned}
\end{equation}

Since $w(v_1)\le\cdots\le w(v_r)$ and
$\sum_{1\le i<j\le r}z'(v_iv_j)\ge1$, we have
\[
\sum_{1\le i<j\le r}z'(v_iv_j)w(v_i)w(v_j)
\ge w(v_1)w(v_2).
\]
By Cauchy--Schwarz,
\[
\begin{aligned}
&(w(v_1)+w(v_2))^2+\sum_{i=3}^r w(v_i)^2\ge
\frac1{r-1}\left(\sum_{i=1}^r w(v_i)\right)^2
=\frac1{r-1}.
\end{aligned}
\]
Consequently,
\begin{equation}\label{eq:fs-clique-bound}
\begin{aligned}
&\sum_{1\le i<j\le r}
(1-z'(v_iv_j))w(v_i)w(v_j)\le
\sum_{1\le i<j\le r}w(v_i)w(v_j)
-w(v_1)w(v_2)\\
&\quad=
\frac12\left(
1-(w(v_1)+w(v_2))^2
-\sum_{i=3}^r w(v_i)^2
\right)\le\frac{r-2}{2(r-1)}.
\end{aligned}
\end{equation}
By \eqref{eq:fs-small-error}, \eqref{eq:fs-weight-comparison},
and the facts that $\varepsilon>0$ and
$0\le z(e),z'(v_iv_j)\le1$, we have
\[
\begin{aligned}
\frac{r-2}{2(r-1)}-\varepsilon^2
&\le\frac1{n^2}\sum_{e\in E(G)}(1-z(e))\le\frac1{n^2}\sum_{e\in E(G)}
(1-z(e))(1+\varepsilon z(e))\\
&\le\sum_{1\le i<j\le r}
(1-z'(v_iv_j))(1+\varepsilon z'(v_iv_j))
w(v_i)w(v_j)\\
&\le(1+\varepsilon)
\sum_{1\le i<j\le r}
(1-z'(v_iv_j))w(v_i)w(v_j).
\end{aligned}
\]
Dividing by $1+\varepsilon>0$ and using
\eqref{eq:fs-clique-bound}, we obtain
\begin{equation}\label{eq:fs-small-defect}
\begin{aligned}
0
&\le\frac{r-2}{2(r-1)}
-\sum_{1\le i<j\le r}
(1-z'(v_iv_j))w(v_i)w(v_j)\\
&\le\frac{r-2}{2(r-1)}- \frac{\frac{r-2}{2(r-1)}
-\varepsilon^2}{1+\varepsilon}=
\frac{\frac{r-2}{2(r-1)}\varepsilon
+\varepsilon^2}{1+\varepsilon}
\le\varepsilon.
\end{aligned}
\end{equation}
The last inequality follows from
$\frac{r-2}{2(r-1)}\le1$ and $\varepsilon>0$.
Expanding the difference in \eqref{eq:fs-small-defect},
we obtain
\[
\begin{aligned}
&\frac{r-2}{2(r-1)}
-\sum_{1\le i<j\le r}
(1-z'(v_iv_j))w(v_i)w(v_j)=\frac12\left[
\left(w(v_1)+w(v_2)-\frac1{r-1}\right)^2
+\sum_{i=3}^r
\left(w(v_i)-\frac1{r-1}\right)^2
\right]\\
&\quad+w(v_1)w(v_2)
\left(\sum_{1\le i<j\le r}z'(v_iv_j)-1\right)+
\sum_{\substack{1\le i<j\le r\\(i,j)\ne(1,2)}}
z'(v_iv_j)
\bigl(w(v_i)w(v_j)-w(v_1)w(v_2)\bigr).
\end{aligned}
\]
Every term on the right is nonnegative, since
$w(v_1)\le\cdots\le w(v_r)$ and
$\sum_{1\le i<j\le r}z'(v_iv_j)\ge1$.
Thus
\[
\left|w(v_1)+w(v_2)-\frac1{r-1}\right|
\le\sqrt{2\varepsilon},
\qquad
\left|w(v_3)-\frac1{r-1}\right|
\le\sqrt{2\varepsilon},
\]
and
\begin{equation}\label{eq:fs-defect-bound}
\begin{aligned}
&\sum_{\substack{1\le i<j\le r\\(i,j)\ne(1,2)}}
z'(v_iv_j)
\bigl(w(v_i)w(v_j)-w(v_1)w(v_2)\bigr)\le\varepsilon.
\end{aligned}
\end{equation}

The preceding bounds imply
$w(v_3)-w(v_2)\ge w(v_1)-2\sqrt{2\varepsilon}$.
Since the weights are ordered, for every
$1\le i<j\le r$ with $(i,j)\ne(1,2)$, we have
\[
\begin{aligned}
w(v_i)w(v_j)-w(v_1)w(v_2)
&\ge w(v_1)(w(v_3)-w(v_2))\ge w(v_1)^2
-2\sqrt{2\varepsilon}\,w(v_1)\\&= \frac{1}{2} w(v_1)^2 + \frac{1}{2}\left( w(v_1)^2 - 4\sqrt{2\varepsilon}\, w(v_1) \right) \\
&= \frac{1}{2} w(v_1)^2 + \frac{1}{2}\left[ \left( w(v_1) - 2\sqrt{2\varepsilon} \right)^2 - 8\varepsilon \right] \\
&\ge \frac{1}{2} w(v_1)^2 - 4\varepsilon.
\end{aligned}
\]

Together with \eqref{eq:fs-defect-bound} and
$0\le z'(v_iv_j)\le1$, this gives
\[
\begin{aligned}
\frac12w(v_1)^2
\sum_{\substack{1\le i<j\le r\\(i,j)\ne(1,2)}}
z'(v_iv_j)
&\le\varepsilon+
4\varepsilon
\sum_{\substack{1\le i<j\le r\\(i,j)\ne(1,2)}}
z'(v_iv_j)\le\left(1+4\binom r2\right)\varepsilon.
\end{aligned}
\]
Since $0\le z'(v_iv_j)\le1$, we also have
\[
\begin{aligned}
&\left(
w(v_1)
\sum_{\substack{1\le i<j\le r\\(i,j)\ne(1,2)}}
z'(v_iv_j)
\right)^2\le
2\binom r2\left(1+4\binom r2\right)\varepsilon.
\end{aligned}
\]
Taking square roots, we obtain
\begin{equation}\label{eq:fs-small-weight}
\begin{aligned}
&w(v_1)
\sum_{\substack{1\le i<j\le r\\(i,j)\ne(1,2)}}
z'(v_iv_j)\le
\sqrt{2\binom r2
\left(1+4\binom r2\right)\varepsilon}\le4\binom r2\sqrt{\varepsilon}.
\end{aligned}
\end{equation}
The last inequality follows from $r\ge3$.
Since $0\le z'(v_iv_j)\le1$ for all
$1\le i<j\le r$, we have
$
z'(v_iv_j)(1-z'(v_iv_j))
\le\min\{z'(v_iv_j),1-z'(v_iv_j)\}.
$
Moreover, $\sum_{1\le i<j\le r}z'(v_iv_j)\ge1$
implies
$
1-z'(v_1v_2)\le
\sum_{\substack{1\le i<j\le r\\(i,j)\ne(1,2)}}
z'(v_iv_j).
$
Consequently,
\[
\begin{aligned}
&\sum_{1\le i<j\le r}
z'(v_iv_j)(1-z'(v_iv_j))w(v_i)w(v_j)\\
&\le w(v_1)w(v_2)(1-z'(v_1v_2))
+\sum_{\substack{1\le i<j\le r\\(i,j)\ne(1,2)}}
z'(v_iv_j)w(v_i)w(v_j)\\
&\le
2w(v_1)w(v_2)
\sum_{\substack{1\le i<j\le r\\(i,j)\ne(1,2)}}
z'(v_iv_j)+
\sum_{\substack{1\le i<j\le r\\(i,j)\ne(1,2)}}
z'(v_iv_j)
\bigl(w(v_i)w(v_j)-w(v_1)w(v_2)\bigr)\le8\binom r2\sqrt{\varepsilon}+\varepsilon.
\end{aligned}
\]
The last inequality follows from
\eqref{eq:fs-defect-bound},
\eqref{eq:fs-small-weight}, and $w(v_2)\le1$.

Finally, by \eqref{eq:fs-weight-comparison},
\eqref{eq:fs-clique-bound}, and the preceding estimate,
we have
\[
\begin{aligned}
&\frac1{n^2}\sum_{e\in E(G)}
(1-z(e))(1+\varepsilon z(e))
\\&\le
\sum_{1\le i<j\le r}
(1-z'(v_iv_j))w(v_i)w(v_j)+
\varepsilon\sum_{1\le i<j\le r}
z'(v_iv_j)(1-z'(v_iv_j))w(v_i)w(v_j)\\
&\le
\frac{r-2}{2(r-1)}
+\varepsilon\left(
8\binom r2\sqrt{\varepsilon}+\varepsilon
\right).
\end{aligned}
\]
Combining this with \eqref{eq:fs-small-error}, we obtain
\[
\begin{aligned}
\frac{\varepsilon}{n^2}
\sum_{e\in E(G)}z(e)(1-z(e))
&=\frac1{n^2}\sum_{e\in E(G)}
(1-z(e))(1+\varepsilon z(e))
-\frac1{n^2}\sum_{e\in E(G)}(1-z(e))\\
&\le\frac{r-2}{2(r-1)}
+\varepsilon\left(
8\binom r2\sqrt{\varepsilon}+\varepsilon
\right)
-\left(\frac{r-2}{2(r-1)}-\varepsilon^2\right)\\
&=8\binom r2\varepsilon^{3/2}
+2\varepsilon^2.
\end{aligned}
\]
Dividing by $\varepsilon>0$ gives
\[
0\le\frac1{n^2}\sum_{e\in E(G)}z(e)(1-z(e))
\le8\binom r2\sqrt{\varepsilon}
+2\varepsilon.
\]
Since $\varepsilon>0$ is arbitrary, the claim follows.
\end{proof}

Let $H$ be the spanning subgraph of $G$ with
$E(H)=\{e\in E(G):z(e)<\binom{r}{2}^{-1}\}$.
Then $H$ is $K_r$-free, since otherwise the sum of
$z(e)$ over the edges of a copy of $K_r$ in $H$
would be less than $1$, contradicting the fractional
$K_r$-cover condition.

By the definition of $H$, every $e\in E(G)\setminus E(H)$
satisfies $z(e)\ge\binom{r}{2}^{-1}$, and hence
$1-z(e)\le\binom{r}{2}z(e)(1-z(e))$.
Therefore, by \eqref{eq:fs-weight-total} and
Claim~\ref{clm:fs-rounding-error}, we obtain
\[
\begin{aligned}
e(H)
&\ge\sum_{e\in E(H)}(1-z(e))=\sum_{e\in E(G)}(1-z(e))
  -\sum_{e\in E(G)\setminus E(H)}(1-z(e))\\
&\ge\sum_{e\in E(G)}(1-z(e))
  -\binom{r}{2}
   \sum_{e\in E(G)\setminus E(H)}z(e)(1-z(e))\\
&\ge\sum_{e\in E(G)}(1-z(e))
  -\binom{r}{2}
   \sum_{e\in E(G)}z(e)(1-z(e))\\
&\ge t_{r-1}(n)-o(n^2).
\end{aligned}
\]
Since $H$ is $K_r$-free, Tur\'an's theorem gives
$e(H)\le t_{r-1}(n)$. Consequently,
$e(H)=t_{r-1}(n)-o(n^2)$.

By Theorem~\ref{thm:furedi}, deleting at most
$t_{r-1}(n)-e(H)$ edges makes $H$ $(r-1)$-partite.
Thus $H$ has an $(r-1)$-partite spanning subgraph
with at least $2e(H)-t_{r-1}(n)$ edges.
For $n\ge r-1$, its vertices can be partitioned
into $r-1$ nonempty independent sets, splitting
independent sets if necessary.

This partition has at least $2e(H)-t_{r-1}(n)$
crossing edges in $G$. Since
$\mathcal P=(V_1,\ldots,V_{r-1})$ maximizes
the number of crossing edges of $G$, we have
$
\sum_{1\le i<j\le r-1}e_G(V_i,V_j)
\ge2e(H)-t_{r-1}(n).
$
By the definitions of $M$ and $D$,
\[
\begin{aligned}
M(\mathcal P)+D(\mathcal P)=t_{r-1}(n)
-\sum_{1\le i<j\le r-1}e_G(V_i,V_j)\le2\bigl(t_{r-1}(n)-e(H)\bigr)
=o(n^2),
\end{aligned}
\]
as required.
\end{proof}

We recall Hall's marriage theorem in the following form.
Let $H$ be a bipartite graph with parts $X$ and $Y$.
Then $H$ has a matching covering every vertex of $X$
if and only if
\[
|N_H(S)|\ge |S|
\qquad\text{for every }S\subseteq X.
\]

\begin{lem}\label{lem:transversal-hall}
Let $r\ge3$ be an integer and $d\ge0$ a real number.
Let $H$ be a graph with an $(r-1)$-partition
$\mathcal P=(V_1,\ldots,V_{r-1})$.
Suppose that $|V_j\setminus N_H(v)|\le d$
for all distinct $i,j\in[r-1]$ and every $v\in V_i$.
If $|V_j|\ge\max\{|V_1|,2(r-2)d\}$
for every $j\in[r-1]\setminus\{1\}$, then $H$ contains
$|V_1|$ pairwise vertex-disjoint copies of $K_{r-1}$,
each meeting every part of $\mathcal P$ in exactly
one vertex.
\end{lem}

\begin{proof}
We construct the required cliques by successively
extending the vertices of $V_1$ through the parts
$V_2,\ldots,V_{r-1}$.

Suppose that, for some $\ell\in[r-2]$, we have
a family $\mathcal K_\ell$ of $|V_1|$ pairwise
vertex-disjoint copies of $K_\ell$, each containing
exactly one vertex from each of $V_1,\ldots,V_\ell$.
For $\ell=1$, take
$\mathcal K_1=\{H[\{v\}]:v\in V_1\}$.
Let $X$ be a set of vertices disjoint from $V(H)$,
in one-to-one correspondence with the cliques
in $\mathcal K_\ell$. For each $C\in\mathcal K_\ell$,
denote the corresponding vertex of $X$ by $x_C$.
Thus $X=\{x_C:C\in\mathcal K_\ell\}$ and
$|X|=|\mathcal K_\ell|=|V_1|$.
Define an auxiliary bipartite graph $Q$ with
parts $X$ and $Y=V_{\ell+1}$, where, for
$C\in\mathcal K_\ell$ and $v\in Y$,
$x_Cv\in E(Q)$ if and only if
$V(C)\subseteq N_H(v)$.

By the definition of $Q$, every $x_C\in X$ satisfies
$|Y\setminus N_Q(x_C)|
\le\sum_{u\in V(C)}|V_{\ell+1}\setminus N_H(u)|
\le\ell d$.
For every $v\in Y$, each clique $C\in\mathcal K_\ell$
with $x_C\notin N_Q(v)$ contains a vertex
not adjacent to $v$ in $H$.
Since the cliques in $\mathcal K_\ell$ are pairwise
vertex-disjoint, these vertices are distinct for
distinct cliques. Hence
$|X\setminus N_Q(v)|
\le\sum_{i\in[\ell]}|V_i\setminus N_H(v)|
\le\ell d$.

We now verify Hall's condition.
Let $S\subseteq X$. The case $S=\varnothing$
is immediate, so assume that $S\ne\varnothing$.
If $|S|\le |Y|-\ell d$, then, for any $x_C\in S$,
$|N_Q(S)|\ge |N_Q(x_C)|
\ge |Y|-\ell d\ge |S|$.
Otherwise, since $|Y|\ge2(r-2)d\ge2\ell d$,
we have $|S|>|Y|-\ell d\ge\ell d$.
As each vertex of $Y$ has at most $\ell d$
nonneighbors in $X$, it has a neighbor in $S$.
Thus $N_Q(S)=Y$, and
$|N_Q(S)|=|Y|\ge|V_1|=|X|\ge|S|$.

By Hall's theorem, $Q$ has a matching covering $X$.
For each $C\in\mathcal K_\ell$, let $y_C\in Y$
be the vertex matched to $x_C$.
By the definition of $Q$,
$H[V(C)\cup\{y_C\}]\cong K_{\ell+1}$.
Since the vertices $y_C$ are distinct and lie
outside $\bigcup_{i\in[\ell]}V_i$, the family
$\mathcal K_{\ell+1}
=\{H[V(C)\cup\{y_C\}]:C\in\mathcal K_\ell\}$
consists of $|V_1|$ pairwise vertex-disjoint
copies of $K_{\ell+1}$, each containing exactly
one vertex from each of $V_1,\ldots,V_{\ell+1}$.
Repeating this construction for
$\ell=1,\ldots,r-2$ gives the required
copies of $K_{r-1}$.
\end{proof}

Let $\mathcal P=(V_1,\ldots,V_{r-1})$ be an
$(r-1)$-partition of $G$.
For an internal edge $xy\in E(G[V_i])$,
where $1\le i\le r-1$, a \emph{transversal extension}
of $xy$ with respect to $\mathcal P$ is a copy of $K_r$
containing $x,y$ and exactly one vertex from each
part $V_j$ with $j\ne i$.

\begin{lem}\label{lem:stars}
Let $G$ be an $n$-vertex graph, and let
$\mathcal P=(V_1,\ldots,V_{r-1})$ be an
$(r-1)$-partition of $G$ maximizing the number
of crossing edges.
Suppose that $|V_i|\ge n/(2(r-1))$
for every $i\in[r-1]$.
Let $B\subseteq V(G)$ satisfy
\[
6|B|^2\le M
\qquad\text{and}\qquad
\mu(v)+2|B|\le\frac{n}{50r^2}
\quad\text{for every }v\in V(G)\setminus B.
\]
Then there exists a family of pairwise
edge-disjoint transversal extensions of
internal edges incident with $B$, each
containing exactly one vertex of $B$,
that covers all but at most $49M/100$
internal edges incident with $B$.
Moreover, every vertex outside $B$ belongs
to at most $|B|$ members of this family.
\end{lem}

\begin{proof}
If $B=\varnothing$, the empty family has the
required properties. Hence assume that
$B\ne\varnothing$.
Since $6|B|^2\le M<n^2/2$, we have
$B\ne V(G)$.
Applying the second hypothesis to any vertex
outside $B$, and using $\mu(v)\ge0$, gives
$|B|\le n/(100r^2)$.

Write $B=\{v_1,\ldots,v_{|B|}\}$ and put
$\mathcal F_0=\varnothing$.
We construct families $\mathcal F_t$,
$0\le t\le |B|$, recursively.
Each $\mathcal F_t$ will consist of pairwise
edge-disjoint transversal extensions, each
meeting $B$ in exactly one of the first $t$
vertices of this ordering.
We also require that every vertex outside $B$
belong to at most $t$ members of $\mathcal F_t$.
These conditions hold for $\mathcal F_0$.
Fix $t\in[|B|]$ and suppose that
$\mathcal F_{t-1}$ has been constructed.
Let $H_t$ be the graph on $V(G)\setminus B$
consisting of the crossing edges unused by
$\mathcal F_{t-1}$.
Since every member of $\mathcal F_{t-1}$
meets $B$ in one of $v_1,\ldots,v_{t-1}$,
none contains $v_t$. Thus no edge incident
with $v_t$ has been used.

By the induction hypothesis, every vertex
outside $B$ belongs to at most $t-1$ members
of $\mathcal F_{t-1}$.
Each transversal extension uses at most two
crossing edges from a vertex into any other part.
Hence, for all distinct $i,j\in[r-1]$
and every $x\in V_i\setminus B$,
\begin{equation}\label{eq:stars-unused-crossing}
|(V_j\setminus B)\setminus N_{H_t}(x)|
\le\mu(x)+2(t-1)
\le\mu(x)+2|B|
\le\frac{n}{50r^2}.
\end{equation}
Let $i\in[r-1]$ satisfy $v_t\in V_i$.
For each $\ell\in[r-1]\setminus\{i\}$,
put $S_\ell=N_G(v_t)\cap(V_\ell\setminus B)$.
Since $\mathcal P$ maximizes the number of
crossing edges, \eqref{eq:local-optimal} holds.
Also, $d_G(v_t,V_\ell)\le |S_\ell|+|B|$
for every $\ell\in[r-1]\setminus\{i\}$.
Therefore,
\begin{equation}\label{eq:stars-local-degree}
\begin{aligned}
d_G(v_t,V_i\setminus B)\le d_G(v_t,V_i)
\le \min_{\ell\in[r-1]\setminus\{i\}}
d_G(v_t,V_\ell)\le \min_{\ell\in[r-1]\setminus\{i\}}
|S_\ell|+|B|.
\end{aligned}
\end{equation}
We will choose a possibly empty family $\mathcal Q_t$
of pairwise vertex-disjoint copies of $K_{r-1}$ in $H_t$.
Each copy will meet every part of $\mathcal P$
in exactly one vertex, and all its vertices will be
adjacent to $v_t$ in $G$.
Together with $v_t$, each copy forms a transversal
extension of an internal edge incident with $v_t$.
We consider two cases.

Suppose first that $|S_j|<(r-2)n/(25r^2)$
for some $j\in[r-1]\setminus\{i\}$.
By \eqref{eq:stars-local-degree}
and $|B|\le n/(100r^2)$,
\[
d_G(v_t,V_i\setminus B)
\le |S_j|+|B|
<\frac{(4r-7)n}{100r^2}.
\]
Since $v_t$ has at most $|S_j|+|B|$
neighbours in $V_j$, the definition of
$\mu(v_t)$ gives
\[
\mu(v_t)\ge |V_j|-|S_j|-|B|
\ge\frac{n}{2(r-1)}
-\frac{(4r-7)n}{100r^2}
\ge\frac{n}{3(r-1)},
\]
where we used $|V_j|\ge n/(2(r-1))$
and $r\ge3$.
Combining the bounds
$d_G(v_t,V_i\setminus B)<(4r-7)n/(100r^2)$
and $\mu(v_t)\ge n/(3(r-1))$, we obtain
\[
d_G(v_t,V_i\setminus B)
\le\frac{3(r-1)(4r-7)}{100r^2}\mu(v_t).
\]
We take $\mathcal Q_t=\varnothing$ in this case,
so no extension is added at step $t$.
All internal edges from $v_t$ to $V_i\setminus B$
remain uncovered, and their number is at most
$3(r-1)(4r-7)\mu(v_t)/(100r^2)$.

Now suppose that $|S_\ell|\ge(r-2)n/(25r^2)$
for every $\ell\in[r-1]\setminus\{i\}$.
Choose $A\subseteq N_G(v_t)\cap(V_i\setminus B)$
with
\[
|A|=\min\left\{
d_G(v_t,V_i\setminus B),
\min_{\ell\in[r-1]\setminus\{i\}}|S_\ell|
\right\}.
\]
If $A=\varnothing$, set $\mathcal Q_t=\varnothing$.
Otherwise, consider the subgraph of $H_t$ induced by
$A\cup\bigcup_{\ell\in[r-1]\setminus\{i\}}S_\ell$.
We apply Lemma~\ref{lem:transversal-hall}
to this induced subgraph with $d=n/(50r^2)$
and the $(r-1)$-partition
\[
(A,S_1,\ldots,S_{i-1},S_{i+1},\ldots,S_{r-1}),
\]
taking $A$ as the first part.

Each of these parts is contained in a distinct
part of $\mathcal P$ and is disjoint from $B$.
Thus \eqref{eq:stars-unused-crossing} implies
that every vertex has at most $d$ nonneighbours
in each other part of this induced subgraph.
Moreover, the choice of $A$ and the assumption
on the sets $S_\ell$ give
\[
|S_\ell|\ge
\max\left\{|A|,\frac{(r-2)n}{25r^2}\right\}
=\max\{|A|,2(r-2)d\}
\]
for every $\ell\in[r-1]\setminus\{i\}$.
Lemma~\ref{lem:transversal-hall} therefore gives
a family $\mathcal Q_t$ of $|A|$ pairwise
vertex-disjoint copies of $K_{r-1}$,
each containing exactly one vertex of $A$
and one vertex of each $S_\ell$,
$\ell\in[r-1]\setminus\{i\}$.
All vertices of these cliques lie in
$N_G(v_t)\setminus B$.

Whether or not $A$ is empty, the family
$\mathcal Q_t$ contains $|A|$ cliques
and meets every vertex of $A$ exactly once.
The corresponding extensions therefore cover
every edge from $v_t$ to $A$.
By the definition of $|A|$ and
\eqref{eq:stars-local-degree}, the number
of uncovered internal edges from $v_t$
to $V_i\setminus B$ in this second case is
\[
\begin{aligned}
d_G(v_t,V_i\setminus B)-|A|=\max\left\{
0,\,
d_G(v_t,V_i\setminus B)
-\min_{\ell\in[r-1]\setminus\{i\}}|S_\ell|
\right\}\le |B|.
\end{aligned}
\]

In either case, every member of $\mathcal Q_t$
is a copy of $K_{r-1}$ in
$H_t[N_G(v_t)\setminus B]$ meeting each part
of $\mathcal P$ in exactly one vertex.
Define
\[
\mathcal F_t=\mathcal F_{t-1}\cup
\bigl\{G[V(Q)\cup\{v_t\}]:Q\in\mathcal Q_t\bigr\}.
\]
Each new member is a transversal extension
whose unique internal edge is incident with $v_t$,
and whose unique vertex in $B$ is $v_t$.
Its edges incident with $v_t$ are unused,
since no member of $\mathcal F_{t-1}$ contains $v_t$.
All its other edges lie in $H_t$ and are
therefore also unused.
Since the cliques in $\mathcal Q_t$ are
pairwise vertex-disjoint, distinct new
extensions intersect only at $v_t$.
Thus $\mathcal F_t$ is pairwise edge-disjoint.
Every vertex outside $B$ belongs to at most
one new extension and hence to at most $t$
members of $\mathcal F_t$.
This completes the inductive step.

It remains to count the uncovered internal
edges incident with $B$.
There are at most $\binom{|B|}{2}$ such edges
with both endpoints in $B$.
Every internal edge with exactly one endpoint
in $B$ is counted at that endpoint.
For a vertex $v_t$ processed in the first case,
at most
$3(r-1)(4r-7)\mu(v_t)/(100r^2)$
such edges remain uncovered.
For a vertex processed in the second case,
at most $|B|$ such edges remain uncovered.
There are at most $|B|$ steps of the second
type, and $\mu(v_t)\ge0$ for every $t$.
Hence the total number of uncovered internal
edges incident with $B$ is at most
\[
\begin{aligned}
\binom{|B|}{2}+|B|^2
+\frac{3(r-1)(4r-7)}{100r^2}
 \sum_{t=1}^{|B|}\mu(v_t)
&\le\frac32|B|^2
+\frac{6(r-1)(4r-7)}{100r^2}M\\
&\le\frac14M+\frac{6}{25}M
=\frac{49}{100}M.
\end{aligned}
\]
The first inequality follows from
$\binom{|B|}{2}\le |B|^2/2$ and
$\sum_{t=1}^{|B|}\mu(v_t)
\le\sum_{v\in V(G)}\mu(v)=2M$,
where the equality is \eqref{eq:sum-deficiency}.
The second follows from $6|B|^2\le M$
and $(r-1)(4r-7)\le4r^2$.

Finally, $\mathcal F_{|B|}$ is pairwise
edge-disjoint, and every vertex outside $B$
belongs to at most $|B|$ members, since
at most one extension containing that vertex
was added at each step.
This proves the lemma.
\end{proof}

\begin{lem}\label{lem:extension}
Fix an integer $r\ge3$ and a real number
$0<\sigma\le1/4$.
For all sufficiently large $q$, the following holds.
Let $H$ be a graph with an $(r-1)$-partition
$\mathcal P=(V_1,\ldots,V_{r-1})$ such that
$\bigl||V_i|-q\bigr|\le\sigma q/100$
for every $i\in[r-1]$.
Suppose that
$|V_j\setminus N_H(x)|\le\sigma q/(100(r-1))$
for all distinct $i,j\in[r-1]$ and every $x\in V_i$.
Let $F$ be a subgraph of $H[V_1]$ with
$\Delta(F)\le(1-\sigma)q$.
Then every edge of $F$ has a transversal extension
with respect to $\mathcal P$, and these extensions
can be chosen to be pairwise edge-disjoint.
\end{lem}

\begin{proof}
If $E(F)=\varnothing$, the conclusion is immediate.
Hence assume that $E(F)\ne\varnothing$.
Take $q\ge2(r-1)D_0(r-1,\sigma/4)$, where
$D_0$ is as in Theorem~\ref{thm:kahn}.

For each $e=uv\in E(F)$, put $X_1(e)=\{u,v\}$.
For $j=2,\ldots,r-1$, we successively choose
$x_j(e)\in V_j$ and put
\[
X_j(e)=X_{j-1}(e)\cup\{x_j(e)\}
=\{u,v,x_2(e),\ldots,x_j(e)\}.
\]
For every $j\in[r-1]$, we require
\begin{equation}\label{eq:extension-invariant}
\begin{aligned}
H[X_j(e)]&\cong K_{j+1}
&&\text{for every }e\in E(F),\\
|X_j(e)\cap X_j(f)|&\le1
&&\text{for all distinct }e,f\in E(F).
\end{aligned}
\end{equation}
These conditions hold for $j=1$, since every
edge of $F$ is an edge of $H$ and distinct
edges of $F$ share at most one endpoint.

Fix $j\in[r-1]\setminus\{1\}$ and suppose that
the sets $X_{j-1}(e)$ have been constructed
and satisfy \eqref{eq:extension-invariant}.
Define a hypergraph $\mathcal H_{j-1}$ by
\[
V(\mathcal H_{j-1})=\bigcup_{\ell\in[j-1]}V_\ell,
\qquad
E(\mathcal H_{j-1})
=\{X_{j-1}(e):e\in E(F)\}.
\]
Each hyperedge has size $j$.
Moreover, the hyperedges are distinct, since
$X_{j-1}(uv)\cap V_1=\{u,v\}$ for every $uv\in E(F)$.
By \eqref{eq:extension-invariant}, any two distinct
hyperedges intersect in at most one vertex.
Thus $\mathcal H_{j-1}$ is linear.

We first bound its maximum degree.
For every $x\in V_1$, the hyperedges containing
$x$ correspond exactly to the edges of $F$
incident with $x$. Hence
\[
d_{\mathcal H_{j-1}}(x)
=\bigl|\{uv\in E(F):x\in\{u,v\}\}\bigr|
\le\Delta(F).
\]
Now let $x\in V_\ell$ for some $2\le\ell\le j-1$.
Since $\mathcal H_{j-1}$ is linear, any two
distinct hyperedges containing $x$ intersect
only at $x$.
As $x\notin V_1$, their intersections with
$V_1$ are pairwise disjoint sets of size two.
Therefore,
\[
2d_{\mathcal H_{j-1}}(x)
=
\left|
\bigcup_{\substack{a\in E(\mathcal H_{j-1})\\x\in a}}
(a\cap V_1)
\right|
\le |V_1|.
\]
Consequently,
\begin{equation}\label{eq:extension-degree}
\begin{aligned}
\Delta(\mathcal H_{j-1})\le\max\left\{\Delta(F),\frac{|V_1|}{2}\right\}\le\max\left\{(1-\sigma)q,
              \frac{1+\sigma/100}{2}q\right\}
=(1-\sigma)q,
\end{aligned}
\end{equation}
where the last equality follows from
$0<\sigma\le1/4$.

Use the vertices of $V_j$ as colours.
For each $a\in E(\mathcal H_{j-1})$, define
\[
L_j(a)=V_j\cap\bigcap_{x\in a}N_H(x).
\]
Thus $L_j(a)$ consists of the vertices of $V_j$
adjacent in $H$ to every vertex of $a$.
Since $a$ has $j$ vertices, all outside $V_j$,
the hypotheses on the part sizes and
crossing nonneighbours give
\begin{equation}\label{eq:extension-list-size}
\begin{aligned}
|L_j(a)|\ge |V_j|-\sum_{x\in a}|V_j\setminus N_H(x)|\ge\left(1-\frac{\sigma}{100}\right)q
       -\frac{j\sigma q}{100(r-1)}\ge\left(1-\frac{\sigma}{50}\right)q
>\frac q2,
\end{aligned}
\end{equation}
where we used $j\le r-1$ and $0<\sigma\le1/4$.

We next find a proper $L_j$-edge-colouring
of $\mathcal H_{j-1}$.
Suppose first that
$\Delta(\mathcal H_{j-1})\le q/(2(r-1))$.
By linearity, for every $a\in E(\mathcal H_{j-1})$,
the number of other hyperedges meeting $a$ is
\[
\begin{aligned}
\bigl|\{b\in E(\mathcal H_{j-1})\setminus\{a\}:
                         a\cap b\ne\varnothing\}\bigr|=\sum_{x\in a}
       \bigl(d_{\mathcal H_{j-1}}(x)-1\bigr)\le(r-1)\bigl(\Delta(\mathcal H_{j-1})-1\bigr)
<\frac q2<|L_j(a)|,
\end{aligned}
\]
where the last inequality follows from
\eqref{eq:extension-list-size}.
Order the hyperedges arbitrarily and colour
them successively.
When colouring a hyperedge $a$, fewer than
$|L_j(a)|$ colours are forbidden by previously
coloured hyperedges meeting $a$.
Thus a colour from $L_j(a)$ remains available,
and this gives a proper $L_j$-edge-colouring.

Otherwise, the choice of $q$ gives
$\Delta(\mathcal H_{j-1})>q/(2(r-1))
\ge D_0(r-1,\sigma/4)$.
For every $a\in E(\mathcal H_{j-1})$,
\eqref{eq:extension-degree} and
\eqref{eq:extension-list-size} imply
\[
\begin{aligned}
(1+\sigma/4)\Delta(\mathcal H_{j-1})
\le(1+\sigma/4)(1-\sigma)q=\left(1-\frac{3\sigma}{4}
          -\frac{\sigma^2}{4}\right)q<\left(1-\frac{\sigma}{50}\right)q
\le |L_j(a)|.
\end{aligned}
\]
Since $\mathcal H_{j-1}$ is linear and has
rank $j\le r-1$, Theorem~\ref{thm:kahn},
with $k=r-1$ and $\gamma=\sigma/4$,
gives a proper $L_j$-edge-colouring.

In either case, let
$\varphi_j:E(\mathcal H_{j-1})\to V_j$
be a proper $L_j$-edge-colouring.
For each $e\in E(F)$, set
\[
x_j(e)=\varphi_j(X_{j-1}(e)),
\qquad
X_j(e)=X_{j-1}(e)\cup\{x_j(e)\}.
\]
Since $x_j(e)\in L_j(X_{j-1}(e))$,
the vertex $x_j(e)$ is adjacent in $H$
to every vertex of $X_{j-1}(e)$.
Also, $x_j(e)\in V_j$ and
$X_{j-1}(e)\cap V_j=\varnothing$.
As $H[X_{j-1}(e)]\cong K_j$ by the induction
hypothesis, we obtain $H[X_j(e)]\cong K_{j+1}$.

It remains to check the intersections.
Let $e,f\in E(F)$ be distinct.
If $X_{j-1}(e)\cap X_{j-1}(f)=\varnothing$,
then $X_j(e)$ and $X_j(f)$ can intersect
only at their newly added vertices in $V_j$,
so $|X_j(e)\cap X_j(f)|\le1$.
Otherwise, properness of $\varphi_j$ gives
$x_j(e)\ne x_j(f)$.
Since neither newly added vertex belongs
to either of the preceding sets, we have
\[
|X_j(e)\cap X_j(f)|
=|X_{j-1}(e)\cap X_{j-1}(f)|
\le1.
\]
Thus both conditions in
\eqref{eq:extension-invariant} are preserved,
completing the induction.

For every $e\in E(F)$, the graph
$H[X_{r-1}(e)]$ is a copy of $K_r$ containing
$e$ and exactly one vertex from each of
$V_2,\ldots,V_{r-1}$.
By \eqref{eq:extension-invariant}, any two
of these copies intersect in at most one
vertex and hence share no edge.
They are therefore the required pairwise
edge-disjoint transversal extensions.
\end{proof}

\section{The two packing estimates}

\begin{lem}\label{lem:sparse}
Fix an integer $r\ge3$, and let $n$ be
sufficiently large.
Let $G$ be a $K_{r+1}$-free graph on $n$ vertices,
and let $\mathcal P=(V_1,\ldots,V_{r-1})$ be an
$(r-1)$-partition of $G$ maximizing the number
of crossing edges.
Suppose that $|V_i|\ge n/(2(r-1))$ for all
$i\in[r-1]$, and that both $I$ and $M$ are
at most $n^2/(10^6r^6)$.
Then $\nu_r(G)\ge I-M/2$.
\end{lem}

\begin{proof}
Assume that $n\ge1000r^3$ and put
\[
\begin{aligned}
B=\left\{v\in V(G):\mu(v)\ge\frac{n}{100r^2}\right\},\ U=\bigcup_{i\in[r-1]}
\left\{v\in V_i\setminus B:
 d_G(v,V_i)\ge\frac{n}{25r^2}\right\}.
\end{aligned}
\]
For each $i\in[r-1]$, put
$C_i=V_i\setminus(B\cup U)$, and let
$C=\bigcup_{i\in[r-1]}C_i$.
Thus $V(G)=B\cup U\cup C$, and these sets
are pairwise disjoint.
All transversal extensions below are taken
with respect to $\mathcal P$.

By the definition of $B$ and \eqref{eq:sum-deficiency},
we have
$n|B|/(100r^2)\le\sum_{v\in B}\mu(v)\le2M$.
Similarly, the definition of $U$ gives
$n|U|/(25r^2)
\le\sum_{i\in[r-1]}\sum_{v\in U\cap V_i}d_G(v,V_i)
\le2I$,
since each internal edge contributes twice to
$\sum_{i\in[r-1]}\sum_{v\in V_i}d_G(v,V_i)$.
Dividing by $n/(100r^2)$ and $n/(25r^2)$,
respectively, and using $I,M\le n^2/(10^6r^6)$,
we obtain
\begin{equation}\label{eq:sparse-exceptional-sizes}
|B|\le\frac{200r^2M}{n}\le\frac{n}{5000r^4},
\qquad
|U|\le\frac{50r^2I}{n}\le\frac{n}{20000r^4}.
\end{equation}

In particular, \eqref{eq:sparse-exceptional-sizes}
gives $|B|+|U|\le n/(4000r^4)$ and
$2|B|+4|U|\le3n/(5000r^4)$.
Since $C_i=V_i\setminus(B\cup U)$ and
$|V_i|\ge n/(2(r-1))$ for every $i\in[r-1]$,
it follows that
\begin{equation}\label{eq:sparse-basic-bounds}
\begin{aligned}
|C_i|
&\ge |V_i|-|B|-|U|
\ge\frac{n}{2(r-1)}-\frac{n}{4000r^4}
\ge\frac{n}{3(r-1)},\\
\frac{n}{100r^2}+2|B|+4|U|
&\le\frac{n}{100r^2}+\frac{3n}{5000r^4}
\le\frac{n}{50r^2}.
\end{aligned}
\end{equation}
The last inequalities follow from
$1/(4000r^4)\le1/(6(r-1))$ and
$3/(5000r^4)\le1/(100r^2)$, respectively,
which hold for $r\ge3$.

We shall repeatedly use the following observation.
For a family $\mathcal A$ of pairwise edge-disjoint
transversal extensions, let $G_{\mathcal A}$ be the
spanning subgraph obtained from $G$ by deleting
all edges of its members.
If $z\in V_i$ belongs to at most $p$ members
of $\mathcal A$, then, for every
$\ell\in[r-1]\setminus\{i\}$,
\begin{equation}\label{eq:sparse-residual}
|C_\ell\setminus N_{G_{\mathcal A}}(z)|
\le\mu(z)+2p.
\end{equation}
Indeed, at most $\mu(z)$ vertices of $C_\ell$
are nonadjacent to $z$ in $G$, and each extension
containing $z$ uses at most two edges from $z$
to $C_\ell$.

\begin{claim}\label{clm:sparse-U-internal}
Every edge of $G[U]$ is internal with respect to
$\mathcal P$.
\end{claim}

\begin{proof}
Suppose that $xy\in E(G)$ with
$x\in U\cap V_i$, $y\in U\cap V_j$,
and distinct $i,j\in[r-1]$.
Since $d_G(x,V_i)\ge n/(25r^2)$ and
$\mu(y)<n/(100r^2)$, we have
\[
\begin{aligned}
|N_G(x)\cap N_G(y)\cap(V_i\setminus B)|
\ge d_G(x,V_i)-\mu(y)-|B|>\frac{3n}{100r^2}-|B|>0,
\end{aligned}
\]
where the last inequality follows from
\eqref{eq:sparse-exceptional-sizes}.
Choose $x'$ in this intersection.
As $d_G(y,V_j)\ge n/(25r^2)$ and
$x,x'\notin B$, we similarly obtain
\[
\begin{aligned}
|N_G(y)\cap N_G(x)\cap N_G(x')\cap(V_j\setminus B)|\ge d_G(y,V_j)-\mu(x)-\mu(x')-|B|>\frac{2n}{100r^2}-|B|>0.
\end{aligned}
\]
Choose $y'$ in this intersection.
The four vertices are distinct and form a $K_4$,
all outside $B$.
For $r=3$, the existence of this $K_4$ contradicts
the assumption that $G$ is $K_4$-free.

Now assume that $r\ge4$ and put $Q=\{x,x',y,y'\}$.
While $Q$ does not meet every part of $\mathcal P$,
choose $\ell\in[r-1]\setminus\{i,j\}$ such that
$Q\cap V_\ell=\varnothing$.
Before each choice, $G[Q]$ is complete,
$Q\cap B=\varnothing$, and $|Q|\le r$.
Since every $v\in Q$ lies outside both $B$ and
$V_\ell$, the definitions of $B$ and $\mu$ give
$|C_\ell\setminus N_G(v)|\le\mu(v)<n/(100r^2)$.
Together with \eqref{eq:sparse-basic-bounds},
this implies
\[
\begin{aligned}
\left|C_\ell\cap\bigcap_{v\in Q}N_G(v)\right|\ge |C_\ell|-\sum_{v\in Q}\mu(v)>\frac{n}{3(r-1)}-\frac{|Q|n}{100r^2}\ge\frac{n}{3(r-1)}-\frac{n}{100r}>0.
\end{aligned}
\]
Choose $z\in C_\ell\cap\bigcap_{v\in Q}N_G(v)$
and replace $Q$ by $Q\cup\{z\}$.
This adds one vertex from a new part, preserves
the completeness of $G[Q]$, and keeps $Q$
disjoint from $B$.
After $r-3$ steps, $|Q|=r+1$, contradicting
the assumption that $G$ is $K_{r+1}$-free.
\end{proof}

We construct nested families
$\mathcal F_0\subseteq\mathcal F_1
\subseteq\mathcal F_2\subseteq\mathcal F_3$
of pairwise edge-disjoint transversal extensions.
The four stages handle internal edges incident
with $B$, edges of $G[U]$, internal edges between
$U$ and $C$, and internal edges within $C$, respectively.

\begin{claim}\label{clm:sparse-pack-B}
There is a family $\mathcal F_0$ of pairwise
edge-disjoint transversal extensions of internal
edges incident with $B$, each containing exactly
one vertex of $B$, that covers all but at most
$49M/100$ internal edges incident with $B$.
Every vertex outside $B$ belongs to at most
$|B|$ members of $\mathcal F_0$.
\end{claim}

\begin{proof}
By \eqref{eq:sparse-exceptional-sizes} and
$M\le n^2/(10^6r^6)$,
\[
6|B|^2\le\frac{240000r^4M^2}{n^2}
\le\frac{6}{25r^2}M\le M.
\]
For every $v\notin B$, the definition of $B$
and \eqref{eq:sparse-basic-bounds} give
$\mu(v)+2|B|<n/(100r^2)+2|B|\le n/(50r^2)$.
The partition $\mathcal P$ also has the optimality
and part-size properties required by
Lemma~\ref{lem:stars}.
That lemma therefore gives the required family.
\end{proof}

Fix $\mathcal F_0$ as in Claim~\ref{clm:sparse-pack-B}.

\begin{claim}\label{clm:sparse-pack-U}
There is a family $\mathcal F_1\supseteq\mathcal F_0$
of pairwise edge-disjoint transversal extensions
such that $\mathcal F_1\setminus\mathcal F_0$
consists of one extension of each edge of $G[U]$,
with all additional vertices in $C$.
Every vertex of $C$ belongs to at most $|U|$
members of $\mathcal F_1\setminus\mathcal F_0$.
\end{claim}

\begin{proof}
If $E(G[U])=\varnothing$, take
$\mathcal F_1=\mathcal F_0$.
Otherwise, order the edges of $G[U]$ as
$e_1,\ldots,e_m$, where $m=e(G[U])$.
We construct their extensions $K_1,\ldots,K_m$
in this order, keeping them pairwise edge-disjoint
and edge-disjoint from $\mathcal F_0$.
All additional vertices will lie in $C$,
and every vertex of $C$ will belong to at most
$|U|$ of the chosen extensions.
These conditions hold before the first step.

Fix $t\in[m]$ and suppose that
$K_1,\ldots,K_{t-1}$ have been constructed.
Let $\mathcal A=\mathcal F_0\cup\{K_1,\ldots,K_{t-1}\}$,
and let $H$ be the spanning subgraph of $G$
obtained by deleting all edges belonging to
members of $\mathcal A$. Thus
$
E(H)=E(G)\setminus\bigcup_{K\in\mathcal A}E(K).
$
By Claim~\ref{clm:sparse-U-internal},
we may write $e_t=xy$ with $x,y\in U\cap V_i$
for some $i\in[r-1]$.
Since $U\cap B=\varnothing$ and every member
of $\mathcal F_0$ has its unique internal edge
incident with $B$, no member of $\mathcal F_0$
contains $e_t$.
For each $1\le s<t$, the unique internal edge
of $K_s$ is $e_s\ne e_t$, so $K_s$ does not
contain $e_t$ either.
Therefore, $e_t$ has not been deleted and
belongs to $E(H)$.

A vertex $z\in U$ can belong to $K_s$ only
as an endpoint of $e_s$, since all additional
vertices lie in $C$.
It therefore belongs to at most
$d_{G[U]}(z)\le|U|-1$ of $K_1,\ldots,K_{t-1}$.
Vertices of $C$ belong to at most $|U|$
of these extensions by the induction hypothesis.
Together with Claim~\ref{clm:sparse-pack-B}, this
shows that every vertex outside $B$ belongs
to at most $|B|+|U|$ members of $\mathcal A$.
By \eqref{eq:sparse-residual}, the definition of $B$,
and \eqref{eq:sparse-basic-bounds},
\[
\begin{aligned}
|C_\ell\setminus N_H(z)|\le\mu(z)+2|B|+2|U|<\frac{n}{100r^2}+2|B|+2|U|
\le\frac{n}{50r^2}
\end{aligned}
\]
for all distinct $k,\ell\in[r-1]$ and
$z\in V_k\setminus B$.

For each $\ell\in[r-1]$, let $T_\ell$ consist
of the vertices of $C_\ell$ belonging to exactly
$|U|$ of $K_1,\ldots,K_{t-1}$.
Each $K_s$, $1\le s<t$, contains at most one
vertex of $C_\ell$, since it is a transversal
extension of an edge of $G[U]$.
Counting incidences between $T_\ell$ and
these cliques gives
$|U||T_\ell|\le t-1\le m\le\binom{|U|}{2}$.
Hence $|T_\ell|\le(|U|-1)/2\le|U|/2$.

Starting with $Q=\{x,y\}$, extend the clique
by choosing one vertex from each $C_\ell$
with $\ell\in[r-1]\setminus\{i\}$.
At the step for $C_\ell$, choose a vertex outside
$T_\ell$ adjacent in $H$ to every vertex of $Q$.
Before this choice, $Q\cap(B\cup V_\ell)=\varnothing$
and $|Q|\le r-1$.
Using the bound $n/(50r^2)$ on the number
of nonneighbours in $C_\ell$, together with
\eqref{eq:sparse-exceptional-sizes} and
\eqref{eq:sparse-basic-bounds}, we obtain
\[
\begin{aligned}
\left|(C_\ell\setminus T_\ell)
\cap\bigcap_{z\in Q}N_H(z)\right|&\ge |C_\ell|-|T_\ell|
-\frac{(r-1)n}{50r^2}\ge\frac{n}{3(r-1)}-\frac{n}{40000r^4}
-\frac{(r-1)n}{50r^2}>0.
\end{aligned}
\]
Thus all $r-2$ choices are possible.
Let $K_t=H[Q]$ for the resulting set $Q$.
This is a transversal extension of $e_t$ and
uses no edge of any previously chosen extension.
Each additional vertex avoids its set $T_\ell$
and hence belonged to at most $|U|-1$
of $K_1,\ldots,K_{t-1}$.
Its number of appearances is therefore at most
$|U|$ after this step; all other vertices of $C$
have unchanged counts.
This completes the induction.
Taking $\mathcal F_1=\mathcal F_0\cup\{K_1,\ldots,K_m\}$
proves the claim.
\end{proof}

Fix $\mathcal F_1$ as in Claim~\ref{clm:sparse-pack-U}.

\begin{claim}\label{clm:sparse-pack-UC}
There is a family $\mathcal F_2\supseteq\mathcal F_1$
of pairwise edge-disjoint transversal extensions
covering all but at most $3|B||U|$ internal
edges between $U$ and $C$.
Each member of $\mathcal F_2\setminus\mathcal F_1$
extends such an edge and contains exactly one
vertex of $U$, with all other vertices in $C$.
Every vertex of $C$ belongs to at most $|U|$
members of $\mathcal F_2\setminus\mathcal F_1$.
\end{claim}

\begin{proof}
If $U=\varnothing$, take $\mathcal F_2=\mathcal F_1$.
Otherwise, write $U=\{u_1,\ldots,u_{|U|}\}$
and process these vertices in order.
At the step for $u_t$, we add extensions containing
$u_t$ and otherwise only vertices of $C$,
and require them to be pairwise vertex-disjoint
outside $u_t$.

Fix $t\in[|U|]$ and suppose that the first
$t-1$ steps have been completed.
Let $u=u_t$, and let $i\in[r-1]$ satisfy $u\in V_i$.
Let $\mathcal A$ consist of $\mathcal F_1$
and all extensions added at these earlier steps.
Let $H$ be the spanning subgraph of $G$
obtained by deleting all edges belonging to
members of $\mathcal A$. Thus
$
E(H)=E(G)\setminus\bigcup_{K\in\mathcal A}E(K).
$
Every extension added at an earlier step of
this stage contains exactly one vertex of $U$,
namely $u_s$ for some $1\le s<t$,
so none of these extensions contains $u$.
Moreover, the unique internal edge of each
member of $\mathcal F_1$ either is incident
with $B$ or has both endpoints in $U$.
Since $u\in U\cap V_i$ and
$C_i\cap(B\cup U)=\varnothing$,
no edge from $u$ to $C_i$ belongs to a member
of $\mathcal F_1$.
Consequently, all edges from $u$ to $C_i$
remain in $H$, and
$N_H(u)\cap C_i=N_G(u)\cap C_i$.

For each $\ell\in[r-1]\setminus\{i\}$,
put $S_\ell=N_H(u)\cap C_\ell$.
By Claim~\ref{clm:sparse-U-internal}, $u$ has no
neighbour in $U\cap V_\ell$, so
$d_G(u,C_\ell)\ge d_G(u,V_\ell)-|B|$.
Since $u\notin B$, it belongs to at most $|B|$
members of $\mathcal F_0$, which use at most
$2|B|$ edges from $u$ to $C_\ell$.
The members of $\mathcal F_1\setminus\mathcal F_0$
containing $u$ are precisely the extensions
of the $d_{G[U]}(u)$ edges incident with $u$
in $G[U]$; each uses exactly one edge from
$u$ to $C_\ell$.
No earlier step of the present stage uses an edge
incident with $u$.
Consequently,
\[
\begin{aligned}
|S_\ell|
&\ge d_G(u,C_\ell)-2|B|-d_{G[U]}(u)\ge d_G(u,V_\ell)-3|B|-d_{G[U]}(u)\\
&>|V_\ell|-\frac{n}{100r^2}-3|B|-|U|
\ge\frac{n}{3(r-1)}.
\end{aligned}
\]
Here the strict inequality uses
$\mu(u)<n/(100r^2)$ and $d_{G[U]}(u)\le|U|$;
the last inequality follows from
\eqref{eq:sparse-exceptional-sizes} and
$|V_\ell|\ge n/(2(r-1))$.

All neighbours of $u$ in $U$ lie in $V_i$,
again by Claim~\ref{clm:sparse-U-internal}.
Thus \eqref{eq:local-optimal} and the bound
$|S_\ell|\ge d_G(u,V_\ell)-3|B|-d_{G[U]}(u)$ give
\[
\begin{aligned}
d_G(u,C_i)
&=d_G(u,V_i)-d_{G[U]}(u)-d_G(u,B\cap V_i)\le d_G(u,V_i)-d_{G[U]}(u)\\
&\le d_G(u,V_\ell)-d_{G[U]}(u)
\le |S_\ell|+3|B|
\end{aligned}
\]
for every $\ell\in[r-1]\setminus\{i\}$.
Choose $A\subseteq N_G(u)\cap C_i$ with
\[
|A|=\min\left\{d_G(u,C_i),
\min_{\ell\in[r-1]\setminus\{i\}}|S_\ell|\right\}.
\]
It follows that
$d_G(u,C_i)-|A|\le3|B|$.

Every vertex of $C$ belongs to at most $|B|$
members of $\mathcal F_0$, at most $|U|$
members of $\mathcal F_1\setminus\mathcal F_0$,
and at most $t-1$ extensions from earlier
steps of the present stage.
By \eqref{eq:sparse-residual}, the definition of $B$,
and \eqref{eq:sparse-basic-bounds},
\[
\begin{aligned}
|C_\ell\setminus N_H(x)|\le\mu(x)+2|B|+2|U|+2(t-1)<\frac{n}{100r^2}+2|B|+4|U|
\le\frac{n}{50r^2}
\end{aligned}
\]
for all distinct $k,\ell\in[r-1]$ and $x\in C_k$.

If $A=\varnothing$, add no extension at this step.
Otherwise, apply Lemma~\ref{lem:transversal-hall}
to the subgraph of $H$ induced by $A$ and the
sets $S_\ell$, $\ell\in[r-1]\setminus\{i\}$,
taking $A$ as the first part and $d=n/(50r^2)$.
These sets are nonempty subsets of distinct
parts of $\mathcal P$ and lie in $C$.
The bound $|C_\ell\setminus N_H(x)|\le d$
therefore gives the required nonneighbour condition.
Also, the choice of $A$ and the bound
$|S_\ell|\ge n/(3(r-1))$ imply
\[
|S_\ell|\ge\max\left\{|A|,\frac{n}{3(r-1)}\right\}
\ge\max\{|A|,2(r-2)d\}
\]
for every $\ell\in[r-1]\setminus\{i\}$.
The lemma gives $|A|$ pairwise vertex-disjoint
copies of $K_{r-1}$, each meeting $A$ and
every $S_\ell$ in exactly one vertex.
Since $A\subseteq N_H(u)$ and
$S_\ell\subseteq N_H(u)$, adding $u$ to each
copy gives a transversal extension in $H$.
These extensions intersect only at $u$ and
use no edge of any member of $\mathcal A$.
They cover every edge from $u$ to $A$,
and each vertex of $C$ belongs to at most
one extension added at this step.

In either case, at most
$d_G(u,C_i)-|A|\le3|B|$ internal edges
from $u$ to $C_i$ remain uncovered.
The construction can therefore be continued
through all $|U|$ steps.
Taking $\mathcal F_2$ to consist of $\mathcal F_1$
and all extensions added in this stage proves
the claim.
\end{proof}

Fix $\mathcal F_2$ as in Claim~\ref{clm:sparse-pack-UC}.
Let $R$ be the graph with vertex set $C$ and
edge set $\bigcup_{i\in[r-1]}E(G[C_i])$.
Then $e(R)\le I\le n^2/(10^6r^6)$.
For $z\in C_i$, the definition of $U$ gives
$d_R(z)\le d_G(z,V_i)<n/(25r^2)$,
so $\Delta(R)<n/(25r^2)$.
Every member of $\mathcal F_2$ has its unique
internal edge incident with $B\cup U$;
hence no edge of $R$ has been used.
Moreover, Claims~\ref{clm:sparse-pack-B},
\ref{clm:sparse-pack-U}, and \ref{clm:sparse-pack-UC}
show that every vertex of $C$ belongs to at most
$|B|+2|U|$ members of $\mathcal F_2$.

\begin{claim}\label{clm:sparse-pack-C}
There is a family $\mathcal F_3\supseteq\mathcal F_2$
of pairwise edge-disjoint transversal extensions
such that $\mathcal F_3\setminus\mathcal F_2$
consists of one extension of each edge of $R$,
all with vertices in $C$.
\end{claim}

\begin{proof}
If $E(R)=\varnothing$, take $\mathcal F_3=\mathcal F_2$.
Otherwise, order the edges of $R$ as
$e_t=x_ty_t$, $t\in[m]$, where $m=e(R)$,
and put $a=\lceil n/(1000r^3)\rceil$.
We construct extensions $K_1,\ldots,K_m$ in this
order, all with vertices in $C$, keeping them
pairwise edge-disjoint and edge-disjoint from
$\mathcal F_2$.
After the first $t$ extensions have been chosen,
write
\[
b_t(z)=\bigl|\{s:1\le s\le t,\,
z\in V(K_s)\setminus\{x_s,y_s\}\}\bigr|
\qquad(z\in C).
\]
We require $b_t(z)\le a$ for every $z\in C$.
Thus only appearances as an additional vertex
are counted; appearances as an endpoint of
$e_s$ are not.
Initially, $b_0(z)=0$ for every $z\in C$.

Fix $t\in[m]$ and suppose that
$K_1,\ldots,K_{t-1}$ have been constructed.
Let $\mathcal A=\mathcal F_2\cup\{K_1,\ldots,K_{t-1}\}$,
and let $H$ be the spanning subgraph of $G$
obtained by deleting all edges belonging to
members of $\mathcal A$. Thus
$
E(H)=E(G)\setminus\bigcup_{K\in\mathcal A}E(K).
$
Since $e_t\in E(R)$ and $\mathcal F_2$ uses
no edge of $R$, no member of $\mathcal F_2$
contains $e_t$.
For each $1\le s<t$, the unique internal edge
of $K_s$ is $e_s\ne e_t$, so $K_s$ does not
contain $e_t$ either.
Consequently, $e_t\in E(H)$.

By Claims~\ref{clm:sparse-pack-B},
\ref{clm:sparse-pack-U}, and \ref{clm:sparse-pack-UC},
every $z\in C$ belongs to at most $|B|+2|U|$
members of $\mathcal F_2$.
Among $K_1,\ldots,K_{t-1}$, the vertex $z$
is an endpoint of the edge being extended
in at most $d_R(z)$ cases.
Its remaining appearances are counted by
$b_{t-1}(z)$, which is at most $a$ by the
induction hypothesis.
Hence every $z\in C$ belongs to at most
$|B|+2|U|+d_R(z)+a$ members of $\mathcal A$.

As $n\ge1000r^3$, we have $a\le n/(500r^3)$.
Using \eqref{eq:sparse-residual},
\eqref{eq:sparse-exceptional-sizes}, and the
bounds $\mu(z)<n/(100r^2)$ and
$d_R(z)<n/(25r^2)$, we obtain
\[
\begin{aligned}
|C_\ell\setminus N_H(z)|\le\mu(z)+2|B|+4|U|+2d_R(z)+2a<\left(\frac{9}{100r^2}+\frac{1}{250r^3}
+\frac{3}{5000r^4}\right)n
\le\frac{n}{10r^2}
\end{aligned}
\]
for all distinct $k,\ell\in[r-1]$ and $z\in C_k$,
where the last inequality uses $r\ge3$.

For each $\ell\in[r-1]$, let
$T_\ell=\{z\in C_\ell:b_{t-1}(z)=a\}$.
Each earlier extension $K_s$ contains at most
one vertex of $C_\ell$ other than $x_s,y_s$.
Consequently,
\[
a|T_\ell|
\le\sum_{z\in C_\ell}b_{t-1}(z)
=\sum_{s=1}^{t-1}
|(V(K_s)\setminus\{x_s,y_s\})\cap C_\ell|
\le t-1\le e(R).
\]
Since $a\ge n/(1000r^3)$ and
$e(R)\le n^2/(10^6r^6)$, it follows that
$|T_\ell|\le n/(1000r^3)$.

Let $i\in[r-1]$ satisfy $x_t,y_t\in C_i$.
Starting with $Q=\{x_t,y_t\}$, extend the clique
through the remaining parts, choosing a vertex
from $C_\ell\setminus T_\ell$ adjacent in $H$
to every vertex already in $Q$.
Before each choice, $Q\subseteq C\setminus C_\ell$
and $|Q|\le r-1$.
The bounds on $|T_\ell|$ and the nonneighbours
in $C_\ell$, together with
$|C_\ell|\ge n/(3(r-1))$ from
\eqref{eq:sparse-basic-bounds}, give
\[
\begin{aligned}
\left|(C_\ell\setminus T_\ell)
\cap\bigcap_{z\in Q}N_H(z)\right|
&\ge |C_\ell|-|T_\ell|
-\frac{(r-1)n}{10r^2}\ge\frac{n}{3(r-1)}-\frac{n}{1000r^3}
-\frac{(r-1)n}{10r^2}>0.
\end{aligned}
\]
Thus all $r-2$ choices are possible.
For the resulting set $Q$, let $K_t=H[Q]$.
This is a transversal extension of $e_t$
using only unused edges.
Each newly chosen vertex $z\in C_\ell$
lies outside $T_\ell$, so
$b_t(z)=b_{t-1}(z)+1\le a$.
For all other vertices, including $x_t,y_t$,
the count is unchanged.
Hence the required conditions are preserved.
After all edges have been extended, take
$\mathcal F_3=\mathcal F_2\cup\{K_1,\ldots,K_m\}$.
\end{proof}

Let $\mathcal F_3$ be as in Claim~\ref{clm:sparse-pack-C}.
Every internal edge of $G$ either is incident
with $B$, has both endpoints in $U$, joins
$U$ to $C$, or has both endpoints in some $C_i$.
The four construction claims show that
$\mathcal F_3$ covers all such edges except
at most $49M/100$ of the first type and
$3|B||U|$ of the third type.
By \eqref{eq:sparse-exceptional-sizes},
\[
3|B||U|\le\frac{600r^2|U|}{n}M
\le\frac{3}{100r^2}M.
\]
Each member of $\mathcal F_3$ contains exactly
one internal edge, and these edges are distinct
because the extensions are pairwise edge-disjoint.
Therefore,
\[
\nu_r(G)\ge|\mathcal F_3|
\ge I-\left(\frac{49}{100}+\frac{3}{100r^2}\right)M
\ge I-\frac12M,
\]
where the last inequality follows from $r\ge3$.
\end{proof}

\begin{lem}\label{lem:dense}
Fix an integer $r\ge3$, and let $n$ be sufficiently large.
Let $G$ be an $n$-vertex $K_{r+1}$-free graph, and let
$\mathcal P=(V_1,\ldots,V_{r-1})$ be an $(r-1)$-partition
maximizing the number of crossing edges.
Suppose that, for each $i\in[r-1]$, the size of $V_i$
differs from $n/(r-1)$ by at most $n/(10^{12}r^{12})$,
and that $M\le n^2/(10^{24}r^{24})$ and
$I\ge n^2/(10^6r^6)$.
Then $\nu_r(G)\ge I-M/2$.
\end{lem}

\begin{proof}
Set
\[
B=\left\{v\in V(G):\mu(v)\ge\frac{n}{10^{10}r^{10}}\right\},
\qquad C_i=V_i\setminus B\quad(i\in[r-1]),
\]
and put $C=V(G)\setminus B$.
By the definition of $B$, every $v\in B$ satisfies
$\mu(v)\ge n/(10^{10}r^{10})$.
Thus $n|B|/(10^{10}r^{10})
\le\sum_{v\in B}\mu(v)\le2M$
by \eqref{eq:sum-deficiency}.
Dividing by $n/(10^{10}r^{10})$ and using
$M\le n^2/(10^{24}r^{24})$, we obtain
\begin{equation}\label{eq:dense-exceptional-size}
|B|\le\frac{2\cdot10^{10}r^{10}M}{n}
\le\frac{2n}{10^{14}r^{14}}.
\end{equation}
Since $C_i=V_i\setminus B$, the assumption on the part
sizes and \eqref{eq:dense-exceptional-size} imply that,
for every $i\in[r-1]$,
\begin{equation}\label{eq:dense-C-size}
\biggl||C_i|-\frac{n}{r-1}\biggr|
\le\frac{2n}{10^{12}r^{12}},
\qquad
\frac{n}{2(r-1)}\le |C_i|\le\frac{2n}{r-1}.
\end{equation}

\begin{claim}\label{clm:dense-internal-structure}
After relabelling the parts, $G[C_1]$ is triangle-free,
$e(G[C_1])\ge n^2/(2\cdot10^6r^6)$, and
$E(G[C_j])=\varnothing$ for every
$j\in[r-1]\setminus\{1\}$.
\end{claim}

\begin{proof}
At most $|B|n$ internal edges are incident with $B$.
Thus \eqref{eq:dense-exceptional-size} and the lower
bound on $I$ give
\[
\sum_{i\in[r-1]}e(G[C_i])
\ge I-|B|n\ge\frac{n^2}{2\cdot10^6r^6}.
\]
Relabel the parts so that
$e(G[C_1])\ge n^2/(2\cdot10^6r^6(r-1))$.

For $j\in[r-1]$ and
$Q\subseteq V(G)\setminus(B\cup V_j)$ with
$1\le |Q|\le r$, the definition of $B$ and
\eqref{eq:dense-C-size} give
\begin{equation}\label{eq:dense-common-neighbours}
\begin{aligned}
\left|C_j\cap\bigcap_{v\in Q}N_G(v)\right|\ge |C_j|-\sum_{v\in Q}\mu(v)>\frac{n}{2(r-1)}-\frac{n}{10^{10}r^9}>0.
\end{aligned}
\end{equation}
A triangle in $G[C_1]$ could therefore be extended
to a $K_{r+1}$ by successively adding a common
neighbour from each other part.
Before each choice, the current clique lies outside
$B$, avoids the chosen part, and has at most $r$ vertices,
so \eqref{eq:dense-common-neighbours} applies.
Hence $G[C_1]$ is triangle-free.

Suppose that $xy\in E(G[C_j])$ for some $j\ne1$.
Then $N_G(x)\cap N_G(y)\cap C_1$ is independent.
Indeed, an edge $uv$ in this set would give a $K_4$
on $\{u,v,x,y\}$.
For $r=3$, this contradicts the assumption that $G$
is $K_4$-free; for $r\ge4$, the same extension process
through the remaining $r-3$ parts gives a $K_{r+1}$.

Every edge of $G[C_1]$ therefore has an endpoint in
$C_1\setminus(N_G(x)\cap N_G(y))$, whose size is at most
$\mu(x)+\mu(y)<2n/(10^{10}r^{10})$.
Since $|C_1|\le n$,
\[
\begin{aligned}
e(G[C_1])\le |C_1|(\mu(x)+\mu(y))<\frac{2n^2}{10^{10}r^{10}}
<\frac{n^2}{2\cdot10^6r^6(r-1)},
\end{aligned}
\]
contrary to the choice of $C_1$.
Thus $G[C_j]$ is edgeless for every
$j\in[r-1]\setminus\{1\}$.
Since $e(G[C_1])=\sum_{i\in[r-1]}e(G[C_i])$,
the claimed lower bound follows.
\end{proof}

Keep this labelling of the parts.
Multiplying the two bounds in
\eqref{eq:dense-exceptional-size} gives
$6|B|^2\le24M/(10^4r^4)\le M$.
For every $v\in C$, the definition of $B$ also gives
\[
\mu(v)+2|B|
<\frac{n}{10^{10}r^{10}}+\frac{4n}{10^{14}r^{14}}
\le\frac{n}{8\cdot10^8r^6}
\le\frac{n}{50r^2}.
\]
Together with the maximality of $\mathcal P$ and
$|V_i|\ge|C_i|\ge n/(2(r-1))$, these bounds verify
the hypotheses of Lemma~\ref{lem:stars}.
It gives a family $\mathcal F$ of pairwise
edge-disjoint transversal extensions of internal
edges incident with $B$, leaving at most $49M/100$
such edges uncovered.
Every vertex of $C$ belongs to at most $|B|$
members of $\mathcal F$.

\begin{claim}\label{clm:dense-extension}
There is a family $\mathcal F'$ of transversal
extensions in $G[C]$, one for each edge of $G[C_1]$,
such that $\mathcal F\cup\mathcal F'$ is pairwise
edge-disjoint.
\end{claim}

\begin{proof}
Let $H$ be the spanning subgraph of $G[C]$
obtained by deleting the edges used by $\mathcal F$.
Thus
\[
V(H)=C,
\qquad
E(H)=E(G[C])\setminus\bigcup_{K\in\mathcal F}E(K).
\]
Each member of $\mathcal F$ has its unique internal
edge incident with $B$, so $H[C_1]=G[C_1]$.

Put
$
q=\frac{n}{r-1},\
\sigma=\frac{(r-1)^2}{8\cdot10^6r^6}.
$
Then $0<\sigma\le1/4$.
Choose a vertex $v$ of maximum degree in $G[C_1]$.
Since $G[C_1]$ is triangle-free, its neighbourhood
is independent, so every edge has an endpoint
outside that neighbourhood.
Hence
$e(G[C_1])\le |C_1|(|C_1|-\Delta(G[C_1]))$.
By Claim~\ref{clm:dense-internal-structure} and
\eqref{eq:dense-C-size},
\[
\begin{aligned}
\Delta(G[C_1])\le |C_1|-\frac{e(G[C_1])}{|C_1|}\le q+\frac{2n}{10^{12}r^{12}}
-\frac{(r-1)n}{4\cdot10^6r^6}
\le(1-\sigma)q,
\end{aligned}
\]
where the last inequality uses
$2/(10^{12}r^{12})\le(r-1)/(8\cdot10^6r^6)$.
Also, \eqref{eq:dense-C-size} gives
\[
\bigl||C_i|-q\bigr|
\le\frac{2n}{10^{12}r^{12}}
\le\frac{(r-1)n}{8\cdot10^8r^6}
=\frac{\sigma q}{100}
\qquad(i\in[r-1]).
\]

For distinct $i,j\in[r-1]$ and $x\in C_i$,
at most $|B|$ members of $\mathcal F$ contain $x$,
and each uses at most two edges from $x$ to $C_j$.
Using the bound on $\mu(x)+2|B|$ established when
applying Lemma~\ref{lem:stars}, we obtain
\[
|C_j\setminus N_H(x)|
\le\mu(x)+2|B|
\le\frac{n}{8\cdot10^8r^6}
=\frac{\sigma q}{100(r-1)}.
\]
Since $\sigma$ depends only on $r$, sufficiently
large $n$ also makes $q$ sufficiently large for
Lemma~\ref{lem:extension}.
Applying that lemma to $H$, the partition
$(C_1,\ldots,C_{r-1})$, and the subgraph $G[C_1]$
gives the required extensions.
Their edges lie in $H$, so none is used by
$\mathcal F$.
\end{proof}

Let $\mathcal F'$ be as in
Claim~\ref{clm:dense-extension}.
By Claim~\ref{clm:dense-internal-structure},
every internal edge not incident with $B$ belongs
to $G[C_1]$.
Thus $\mathcal F\cup\mathcal F'$ leaves at most
$49M/100$ internal edges uncovered.
Since these extensions are pairwise edge-disjoint
and each contains exactly one internal edge,
\[
\nu_r(G)\ge|\mathcal F\cup\mathcal F'|
\ge I-\frac{49}{100}M
\ge I-\frac12M,
\]
as required.
\end{proof}

\section{Proof of the main theorem}\label{proof}

\begin{proof}[Proof of Theorem~\ref{thm:main}]
Fix an integer $r\ge3$.
Suppose, for a contradiction, that for arbitrarily large
$n$, there exists an $n$-vertex $K_{r+1}$-free graph $G$
satisfying $\nu_r(G)\le k:=e(G)-t_{r-1}(n)$,
but either $G\notin\mathcal E_r(n)$ or $\nu_r(G)<k$.
Choose an $(r-1)$-partition
$\mathcal P=(V_1,\ldots,V_{r-1})$ maximizing the number
of crossing edges.
By Lemma~\ref{lem:fractional-stability},
we have $M+D=o(n^2)$.

Since $\sum_{i\in[r-1]}|V_i|=n$, the definition of $D$
gives
\[
\begin{aligned}
\frac12\sum_{i\in[r-1]}
\left(|V_i|-\frac{n}{r-1}\right)^2
&=\frac12\sum_{i\in[r-1]}|V_i|^2
-\frac{n^2}{2(r-1)}\\
&=D+\frac{r-2}{2(r-1)}n^2-t_{r-1}(n)=D+O(1)=o(n^2),
\end{aligned}
\]
where we used
$t_{r-1}(n)=\frac{r-2}{2(r-1)}n^2+O(1)$.
Consequently, for all sufficiently large counterexamples,
\[
M\le\frac{n^2}{10^{24}r^{24}},
\qquad
\left||V_i|-\frac{n}{r-1}\right|
\le\frac{n}{10^{12}r^{12}}
\quad(i\in[r-1]).
\]
In particular, every part has size at least
$n/(2(r-1))$.

Apply Lemma~\ref{lem:sparse} if
$I\le n^2/(10^6r^6)$, and Lemma~\ref{lem:dense}
otherwise.
In either case, $\nu_r(G)\ge I-M/2$.
By the definitions of $I,M,D$, we have $k=I-M-D$.
Hence
\[
k\ge\nu_r(G)\ge I-\frac12M
=k+\frac12M+D.
\]
Since $M,D\ge0$, this forces $M=D=0$ and $\nu_r(G)=k$.
Thus $\mathcal P$ is balanced and all crossing edges
are present.

A triangle in some $G[V_i]$, together with one vertex
from each of the other $r-2$ parts, would form a
$K_{r+1}$.
Likewise, two internal edges in distinct parts,
together with one vertex from each of the remaining
$r-3$ parts, would form a $K_{r+1}$.
Since $G$ is $K_{r+1}$-free, all internal edges lie
in one part and form a triangle-free graph.
Thus $G\in\mathcal E_r(n)$.
Therefore, for all sufficiently large $n$, the required
lower bound holds, and equality implies
$G\in\mathcal E_r(n)$.

Conversely, let $G\in\mathcal E_r(n)$, with $n$
sufficiently large.
Choose the balanced partition defining $G$ so that
all added edges lie in $V_1$.
Then $G[V_1]$ is triangle-free and all other parts
are independent.
Every clique therefore contains at most two vertices
from $V_1$ and at most one from each other part.
Thus $G$ is $K_{r+1}$-free, and every copy of $K_r$
contains exactly one internal edge.
It follows that
$\nu_r(G)\le e(G)-t_{r-1}(n)$.
Together with the lower bound already proved, this
gives $\nu_r(G)=e(G)-t_{r-1}(n)$.
\end{proof}

The assumption that $n$ is sufficiently large
in Theorem~\ref{thm:main} is necessary.
The following construction shows that the threshold
for $n$ must grow at least quadratically with $r$.

\begin{pro}\label{pro:quadratic-obstruction}
For every integer $r\ge5$, there exists a
$K_{r+1}$-free graph $G$ on
$n=(r-1)(r-2)+1$ vertices such that
$e(G)=t_{r-1}(n)+2r-6$ but
$\nu_r(G)<2r-6$.
\end{pro}

\begin{proof}
Let $G$ be the complete $r$-partite graph
with parts $V_1,\ldots,V_r$ satisfying
$|V_1|=2$, $|V_2|=r-3$, and
$|V_i|=r-2$ for $3\le i\le r$.
Clearly, $G$ is $K_{r+1}$-free and
\[
n=2+(r-3)+(r-2)^2=(r-1)(r-2)+1.
\]
Combining $V_1$ and $V_2$ gives a balanced
$(r-1)$-partition with part sizes
$r-1,r-2,\ldots,r-2$.
Deleting the $2(r-3)$ edges between
$V_1$ and $V_2$ therefore gives
$T_{r-1}(n)$, so
$e(G)=t_{r-1}(n)+2r-6$.

Suppose, for a contradiction, that $G$ contains
a family $\mathcal K$ of $2r-6$ pairwise
edge-disjoint $K_r$.
Write $V_1=\{a,b\}$, and let
$\mathcal K_a$ and $\mathcal K_b$ consist
of the cliques containing $a$ and $b$,
respectively.
Since every $K_r$ contains exactly one vertex
from each part,
$|\mathcal K_a|+|\mathcal K_b|=2r-6$.

For every $i\in\{2,\ldots,r\}$, the $r-3$
cliques in $\mathcal K_a$ contain $r-3$
distinct vertices of $V_i$, since any two
of them sharing a vertex of $V_i$ would
also share an edge incident with $a$.
The same holds for $\mathcal K_b$.
Hence, by inclusion-exclusion, at least
$2(r-3)-|V_i|$ vertices of $V_i$ each belong
to exactly one clique in $\mathcal K_a$
and exactly one clique in $\mathcal K_b$.
Each such vertex belongs to a unique
$K\in\mathcal K_a$ and a unique
$K'\in\mathcal K_b$.
Moreover, two distinct vertices cannot
belong to the same pair $K,K'$, since
otherwise $K$ and $K'$ would share the edge
joining these vertices, contradicting
edge-disjointness.
Thus the total number of these vertices
over all $V_2,\ldots,V_r$ is at most
$|\mathcal K_a||\mathcal K_b|$.
Consequently,
\[
\begin{aligned}
(r-3)^2=|\mathcal K_a||\mathcal K_b|\ge\sum_{i=2}^r\bigl(2(r-3)-|V_i|\bigr)=(r-3)+(r-2)(r-4)=(r-3)^2+(r-4),
\end{aligned}
\]
a contradiction since $r\ge5$.
Therefore $\nu_r(G)<2r-6$.
\end{proof}

\section{Concluding remarks}

We proved that, for every fixed integer $r\ge3$
and all sufficiently large $n$, every $n$-vertex
$K_{r+1}$-free graph $G$ with $t_{r-1}(n)+k$ edges
satisfies $\nu_r(G)\ge k$.
Moreover, equality holds if and only if
$G\in\mathcal E_r(n)$.

As noted in Remark~1, such a graph satisfies
$\nu_r(G)\ge k+1$ whenever
$\lfloor\lceil n/(r-1)\rceil^2/4\rfloor+1
\le k\le t_r(n)-t_{r-1}(n)$.
For $k$ in this range, it would be interesting
to determine the minimum possible value of
$\nu_r(G)$ among all $n$-vertex $K_{r+1}$-free
graphs with $t_{r-1}(n)+k$ edges, and to
characterize the graphs attaining this minimum.

A further question concerns the order of the least
threshold on $n$ in Theorem~\ref{thm:main}
as a function of $r$.
For $r\ge5$, Proposition~\ref{pro:quadratic-obstruction}
gives the lower bound $(r-1)(r-2)+2$,
but our proof yields no useful explicit upper bound.

\section*{Declaration on the Use of Generative AI}
The authors used ChatGPT 6 Pro to assist in discussing proof strategies, checking proofs, and
improving exposition.

\end{document}